\documentclass[11pt]{article}
\usepackage{amsmath,amssymb,amscd,latexsym,amsthm,mathrsfs}
\usepackage[unicode]{hyperref}
\usepackage{hypbmsec,longtable}
\usepackage{graphicx}
\usepackage{fancyhdr}
\usepackage{topcapt}
\usepackage{subfigure}
\ifx\pdfoutput\undefined\else
\hypersetup{
colorlinks=true,
linkcolor=blue,
citecolor=blue,
urlcolor=blue,
filecolor=blue,
bookmarksnumbered=true,
pdfstartview=FitH,
pdfhighlight=/N
}
\fi

\theoremstyle{remark}

\newcommand{\BE}{\begin{equation}}
\newcommand{\EE}{\end{equation}}

\usepackage{cite}

\renewcommand{\author}[1]{\large\rm #1\\ \bigskip}
\newcommand{\address}[1]{{\normalsize\it #1\\}\bigskip}
\renewcommand{\title}[1]{\bigskip\bigskip\Large\bf #1\bigskip\bigskip\\}

\begin{document}


\vglue .3 cm

\begin{center}

\title{On the growth rate of 1324-avoiding permutations}
\author{Andrew R. Conway\footnote[1]{email:  {\tt andrew1324@greatcactus.org}  }   and Anthony J. Guttmann\footnote[2]{email:
                {\tt tonyg@ms.unimelb.edu.au}}      }

\address{ ARC Centre of Excellence for\\
Mathematics and Statistics of Complex Systems,\\
Department of Mathematics and Statistics,\\
The University of Melbourne, Victoria 3010, Australia}

\end{center}
\setcounter{footnote}{0}

\begin{abstract}
We give an improved algorithm for counting the number of $1324$-avoiding permutations, resulting in 5 further terms of the generating function. We analyse the known coefficients and find compelling evidence that unlike other classical length-4 pattern-avoiding permutations, the generating function in this case does not have an algebraic singularity. Rather, the number of 1324-avoiding permutations of length $n$  behaves as $$B\cdot \mu^n \cdot \mu_1^{n^{\sigma}} \cdot n^g.$$ We estimate $\mu=11.60 \pm 0.01,$ $\sigma=1/2,$ $\mu_1 = 0.0398 \pm 0.0010,$ $g = -1.1 \pm 0.2$ and $B =9.5 \pm 1.0.$

\end{abstract}

\section{Introduction}
Let $\pi$ be a permutation on $[n]$ and $\tau$ be a permutation on $[k].$ Then $\tau$ is said to occur as a {\em pattern} in $\pi$ if for some subsequence  of $\pi$ of length $k$ all the elements of the subsequence occur in the same relative order as do the elements of $\tau.$ For example $1324$ occurs as a pattern in $152364$ as $1526$ and $1536$ as both are in the same relative order as $1324.$ If a permutation $\tau$ does not occur in $\pi,$ then this is said to be a {\em pattern-avoiding permutation,} or PAP.
Let $P(z) = \sum_{n \ge 0} p_nz^n$ be the ordinary generating function (OGF) for the number of permutations $p_n$ of length $n$ avoiding the pattern $1324.$ 
It is well known that, for the classical, length 4 PAPs, the 24 possible patterns  fall into one of three possible classes \cite{ES10}, called Wilf classes. That is to say, there are three distinct OGFs describing all 24 patterns. 

For the sequence $1234$ and its associated patterns, in 1990 Gessel \cite{IG90} showed that the number of length $n >0$ pattern-avoiding permutations is
\begin{equation} \label{eq:1234}
p_n(1234) = \frac{1}{(n+1)^2(n+2)}\sum_{k=0}^n \binom{2k}{k} \binom{n+1}{k+1} \binom{n+2}{k+1}.
\end{equation}
Asymptotically, $$p_n(1234) \sim 2.8 \cdot 9^n \cdot n^{-4},$$

and the generating function $P_{1234}(x) = \sum_n p_n(1234)x^n$ satisfies the linear ODE
\begin{multline}
(9x^5-19x^4+11x^3-x^2)\cdot \frac{d^3P_{1234}(x)}{dx^3} +(72x^4-153x^3+90x^2-9x)\cdot \frac{d^2P_{1234}(x)}{dx^2} + \\
+(126x^3-264x^2+154x-16)\cdot \frac{dP_{1234}(x)}{dx} + (32-72x+36x^2)\cdot P_{1234}(x)=0, \,\,\,\,\,\,\,\,\,\,\,\,  \\
 {\rm with\,\,\, initial\,\,\, conditions\,\,\,}  P_{1234}(0)=1, \,\, P'_{1234}(2)=0, \,\,  P''_{1234}(2)=12.\,\,\,\,\,\,\,\,\,\,\,\,\,\,\,\,\,\,\,\,\,\,\,\,\,\,\,\,\,\,\,\,\,\,\,\,\,\,\,
\end{multline}

For the sequence $1342$ and its associated patterns, in 1997 B\'ona \cite{MB97} showed that the number of length $n >0$ pattern-avoiding permutations is
\BE \label{eq:1342}
p_n(1342) = (-1)^{n-1} \cdot \frac{(7n^2-3n-2)}{2} + 3\sum_{k=0}^n (-1)^{n-i}\cdot 2^{i+1} \cdot \frac{(2i-4)!}{i!(i-2)!} \cdot\binom{n-i+2}{2}.
\EE
The generating function $P_{1342}(x) = \sum_n p_n(1342)x^n$ satisfies the linear ODE
\BE
(8x^2+7x-1)\cdot \frac{d^2P_{1342}(x)}{dx^2}+(28x-8)\cdot \frac{dP_{1342}(x)}{dx} + 12\cdot P_{1342}(x)=0, \,\, P_{1342}(0)=1, \,\, P'_{1342}(0)=1.
\EE
Indeed, it can be exactly solved to give \cite{K11} the simple algebraic expression
$$P_{1342}(x)=\frac{32x}{-8x^2+20x+1-(1-8x)^{3/2}},$$ from which one readily obtains $$p_n(1342) \sim \frac{64}{243 \cdot \sqrt{\pi}} \cdot 8^n \cdot n^{-5/2}.$$

The remaining class, that of 1324-avoiding permutations, remains unsolved \cite{ES12}. Even the growth constant is not accurately known. The best upper bound is 13.73718 \cite{MB14}, due to B\'ona, while Claesson, Jel\'inek and Steingr\'imsson \cite{CJS12} gave an improved bound, $e^{\pi \sqrt{2/3}} \approx 13.00195,$ but subject to the validity of an unproved conjecture. The best published lower bound is 9.47, proved by Albert et al. \cite{AR06} while David~Bevan has an as yet unpublished bound of 9.81. Careful Monte Carlo work by Madras and Liu \cite{ML10} implies that the growth constant lies in the range $[10.71,11.83].$

In this paper we give details of an improved algorithm for the enumeration of such PAPs, with which we obtained five further terms in the OGF beyond the existing longest known sequence, due to Johansson and Nakamura \cite{JN13}, using comparable computing resources. We will refer to their algorithm as the JN algorithm.

We also analyse the sequence of coefficients, and provide compelling numerical evidence that the asymptotic form of the coefficients is more complex than that of the two solved classes just considered. The generating functions of the two solved cases have algebraic singularities, whose coefficients have the asymptotic form $p_n \sim D \cdot \mu^n \cdot n^g.$ Our numerical studies, detailed below, lead us to suggest that the coefficients $p_n(1324) \sim B \cdot \mu^n \cdot \mu_1^{n^\sigma} \cdot n^g,$ where $\sigma = \frac{1}{2}.$ 

Many enumeration problems in algebraic combinatorics have generating functions with algebraic singularities, and hence coefficients with leading asymptotic form 
\BE \label{eq:alg}
a_n \sim A \cdot \mu^n \cdot n^g,
\EE
 where $1/\mu$ is the radius of convergence,  sometimes called the {\em critical point}, $g$ is a critical exponent and $A$ is a critical amplitude. However
a number of solved, and, we believe, unsolved problems that arise in both  algebraic combinatorics and mathematical physics have a more complex singularity structure, with a sub-dominant asymptotic term ${\rm O}(\mu_1^{n^\sigma})$ rather then ${\rm O}(n^g).$ In fact the sub-sub dominant term is of ${\rm O}(n^g)$. That is to say, the dominant asymptotic behaviour of the coefficients $b_n$ of the associated generating function  is
\BE \label{eq:can1}
b_n \sim B \cdot \mu^n \cdot \mu_1^{n^\sigma} \cdot n^{g}.
\EE
Perhaps the best-known example of this sort of behaviour is the number of partitions of the integers -- though in that case the leading exponential growth term $\mu^n$ is absent (or, equivalently, $\mu=1$). Another example is the generating function for the number of fragmented permutations \cite{FS09}, which is $$F(z) =\exp\left ( \frac{z}{1-z} \right ).$$ Then, with $F_n=[z^n]F(z),$ we have\cite{FS09}, p563, $$F_n \sim \frac{e^{2\sqrt{n}}}{2\sqrt{\pi e}n^{3/4}}.$$ These two examples highlight the fact that singularities of quite different analytic structure can give rise to the same asymptotics. The OGF for integer partitions which has radius of convergence $r_c=1,$ has a natural boundary on the unit circle, whereas the generating function for fragmented permutations is D-finite. So if such an asymptotic form (\ref{eq:can1}) is observed, one cannot say much about the underlying singularity structure from the asymptotics alone.

Another example, discussed in \cite{G14}, is the problem of Dyck paths counted not only by length but also by height $h$, which is defined to be the maximal vertical displacement of the Dyck path from its horizontal axis. Let $d_{n,h} $ be the number of Dyck paths of length $2n$ and height $h.$ The generating function is 
$$D(x,y)=\sum_{n,h}d_{n,h}x^{2n}y^h.$$
Then
\BE \label{eq:Dyck}
[x^{2n}]D(x,y)=\sum_{h=1}^n d_{n,h}y^h \sim B\cdot 4^n \cdot \mu_1^{n^{1/3}} \cdot n^{-5/6},
\EE
where both $B(y)$ and $\mu_1(y)$ are known \cite{G14}.

There are also a number of models in the mathematical physics literature that  have this more complex asymptotic structure. In particular, Duplantier and Saleur \cite{DS87} and Duplantier and David \cite{DD88} studied the case of {\em dense} polymers in two dimensions, and found the partition functions had the asymptotic form (\ref{eq:can1}). In \cite{OPB93}, Owczarek, Prellberg and Brak investigated an exactly solvable model of interacting partially-directed self-avoiding walks (IPDSAW), for which the solution had previously been given by Brak, Guttmann and Whittington in \cite{BGW92}. In \cite{OPB93}, Owczarek et al. analysed a 6000 term series expansion for IPDSAWs in the collapsed regime, and estimated $\sigma = 1/2,$ $g = -3/4,$ while $\mu_1$ was estimated to at least 6 digit accuracy. From \cite{BGW92} the value of $\mu$ is exactly known. Subsequently Duplantier \cite{D93} pointed out that $\sigma = 1/2$ is to be expected, not only for IPDSAWs, but also for SAWs in the collapsed regime, for the two-dimensional version of these models. In all the examples we have encountered, $\sigma$ takes the value 1/3, 1/2 or 2/3.

In the next section we give details of the enumeration algorithm. In subsequent sections we analyse the available series coefficients.

\section{Algorithm}

The algorithm used can be considered to be a set of further optimizations on  the JN algorithm.
However a significantly different notation is used as this helps make some of the optimizations
clearer, as well as helping with a memory efficient encoding implementation. We note that Marinov and Rodoi\v{c}i\'c \cite{MR03} have previously given a recursive algorithm for this problem. It is a significantly different algorithm conceptually to the JN algorithm, although one could imagine a variant on that algorithm keeping track of 12 patterns as part of their labels, instead of their set of $l(\pi)$ values. Such an algorithm would probably be similar in performance to this algorithm.

\subsection{Basic Algorithm}

We will start with a very simple (and inefficient) algorithm. Let $f(n,P)$ be the number
of permutations avoiding the pattern 1324 with $n$ numbers remaining and starting with
the prefix $P$ (a sequence of integers). Then the desired series is $f(n,\varnothing)$. Each
value $f(n,P)$ can be expressed as the sum of up to $n$ other terms $f(n-1,P')$ where $P'$ is
$P$ followed by one extra integer. There will be fewer than $n$ terms if $P'$ implies a
1324 pattern, usually due to containing a 132 pattern with a 4 inevitably to eventually follow. 
Using the termination condition $f(0,-)=1$, one could easily
encode a recursive algorithm that would work.

This algorithm will systematically individually enumerate every permutation avoiding 1324,
and its time consumption will be proportional to the answer. Like many such recursive
enumeration algorithms, one can get a much faster algorithm by recognising that there
exist many classes $S = (n,{P_1,P_2,P_3,...})$ such that $f(n,P_i)=f(n,P_j)$ for all $i$ and $j$.
Define $f(S)$=$f(n,P_-), $ $P_-$ means any prefix from the class. Now, modify the algorithm to use $S$ (which will henceforth
be called a signature) in a recursive function $f(S)$. After computing
a value of $f(S)$, store it in some table. When you next need $f(S)$, look it up in the
table. If it is already there, then use the stored value. This can vastly improve
the speed of execution as you will avoid passing through large swathes of the
enumeration tree. It does have a cost of memory. This approach is often called dynamic
programming, memoization, or memorization.

When using such approaches, the definition of the signatures $S$ is paramount. The
more prefixes you can prove to be identical (and thus members of the same signature),
the more efficient the algorithm will be. Indeed, the time and memory will both
be proportional to the number of different signatures. The rest of this subsection deals
with the definition of prefixes.

The signature must contain enough information to allow the algorithm to avoid
1324 patterns. One way to do this is to keep track of all the 13 patterns. Then,
in the recurrence relation, you do not allow any numbers in the middle of
one of these 13 patterns (i.e. a 2) if there are any numbers remaining bigger than
the 3 (i.e. a 4, which would then inevitably follow at some point).
There is no need to keep track of exactly what the numbers are; you just need
to know the total number of numbers between, above, and below each 13 pattern.

One suitable notation to keep track of the 13 patterns is as a partition of
integers (the number of numbers left to go) with a well formed set of brackets
(the 13 patterns). It is also necessary to keep track of the lowest number
so far in the prefix, as all future numbers higher than it will form a 13
pattern with it. They may of course form other 13 patterns, but the one with the
lowest 1 will be the most restrictive, and the others may be ignored. The lowest
number is recorded with a comma. For brevity, the comma may be left out if there is a bracket
immediately following it.

Some typical signatures are shown in table \ref{tab:egSigs}, for some example prefixes and initial $n=20$.

\begin{table}
\caption{Example prefixes and corresponding signatures for $n=20$.}
\label{tab:egSigs}
\begin{center}
\begin{tabular}{lll} \hline \hline

Prefix                   &     Signature    & Explanation \\

\hline \hline

$\varnothing$            &     20           & 20 numbers left to go \\

11                       &     10,9         & ten numbers to the left of the lowest number, 9 to the right \\

11,14                    &     10,[2]6      & a 13 pattern introduced with two numbers between it, and 6 to the right \\

11,14,5                  &     4,5[2]6      & a new lowest number, but no new 13 pairs \\

11,14,5,9                &     4,[3]1[2]6   & a new 13 pair produced \\

11,14,5,9,15             &     4,[[3]1[2]]5 & a new 13 pair produced \\

\hline \hline

\end{tabular}

\end{center}

\end{table}

At this point one can simplify the signature. A signature of the form $a[b[c]]d$ (that is,
with two consecutive closing brackets) has two restrictions.
The outer bracket means that you can't have anything in $b$ or $c$
until $d$ is finished. The inner bracket means that you can't have anything in $c$ until
$d$ is finished. The outer bracket is strictly more restrictive, so the inner one
is redundant and may be removed. So $a[b[c]]d=a[bc]d$. This simplification means that we will never have
two closing brackets in succession.

Adjacent integers not separated by brackets or a comma can be added together. For instance, $2[3[4]]5$, after removing the inner bracket, becomes 
$2[7]5$ rather than needing to record that the inner $7$ was at one point broken unto a $3$ and a $4$.

Any tail bracket at the end of the signature of course can be removed; the brackets
are only restrictive if there are larger numbers possible.

At this point there is an isomorphism to the JN algorithm \cite{JN13}. Their final functional
form (last equation in section 2) is $H^0_n (t; b_1, ... , b_n; k)$. Here $k$ encodes the
position of the comma, and $b_i$ encodes the location of the open bracket corresponding
to a closing bracket at position $i$. Enumerating using these signatures produces an algorithm
basically identical in performance to  the JN algorithm. The further simplifications described
below will improve performance.

Repeated open brackets can also be simplified. Consider a signature of the form $a[[b]c]d$. The
outer bracket means you cannot have anything in $b$ or $c$ until everything in $d$ is done.
The inner bracket means you cannot have anything in $b$ until $c$ is done. This is equivalent
to the restrictions described by the signature $a[b][c]d$. Simplifying signatures
by getting rid of all consecutive open brackets reduces the total number of signatures
significantly, in practice by a factor of roughly 4.

A minor simplification comes from dealing with open brackets at the start of the
signature. $[a]b$ means that everything in $b$ must be dealt with before everything in $a$.
This means that they are decoupled; indeed $f([a]b) = f(a)f(b)$. Factorizing the problem
seems like a big advantage, but only a small proportion of signatures start with brackets;
in practice this reduces the number of signatures by a factor of roughly 2.

An example of all the computations done (in the order that they are finished)
for permutations of length 6 is given in table \ref{tab:exhaustiveList}

\begin{table}

\caption{ All signatures used to compute up to the $n=6$ term. Single lines are used when a higher $n$ is started. }

\label{tab:exhaustiveList}

\begin{center}

\begin{tabular}{lll} \hline \hline

Signature $S$            &  Composed of & $f(S)$ \\

\hline \hline

1 & , & 1 \\

\hline

,1 & , & 1 \\

2 & ,1 + 1 & 2 \\

\hline

,2 & ,1 + ,1 & 2 \\

1,1 & ,1 + 1 & 2 \\

3 & ,2 + 1,1 + 2 & 6 \\

\hline

,3 & ,2 + [1]1 + ,2 & 5 \\

1,2 & ,2 + 1,1 + 1,1 & 6 \\

2,1 & ,2 + 1,1 + 2 & 6 \\

4 & ,3 + 1,2 + 2,1 + 3 & 23 \\

\hline

,4 & ,3 + [1]2 + [2]1 + ,3 & 14 \\

1[1]1 & [1]1 + 1,1 & 3 \\

1,3 & ,3 + 1,2 + 1[1]1 + 1,2 & 20 \\

2,2 & ,3 + 1,2 + 2,1 + 2,1 & 23 \\

3,1 & ,3 + 1,2 + 2,1 + 3 & 23 \\

5 & ,4 + 1,3 + 2,2 + 3,1 + 4 & 103 \\

\hline

,5 & ,4 + [1]3 + [2]2 + [3]1 + ,4 & 42 \\

1[1]2 & [1]2 + 1[1]1 + 1[1]1 & 8 \\

1[2]1 & [2]1 + 1,2 & 8 \\

1,4 & ,4 + 1,3 + 1[1]2 + 1[2]1 + 1,3 & 70 \\

,1[1]1 & [1]1 + ,2 & 3 \\

2[1]1 & ,1[1]1 + 1[1]1 + 2,1 & 12 \\

2,3 & ,4 + 1,3 + 2,2 + 2[1]1 + 2,2 & 92 \\

3,2 & ,4 + 1,3 + 2,2 + 3,1 + 3,1 & 103 \\

4,1 & ,4 + 1,3 + 2,2 + 3,1 + 4 & 103 \\

6 & ,5 + 1,4 + 2,3 + 3,2 + 4,1 + 5 & 513 \\

\hline \hline

\end{tabular}

\end{center}

\end{table}

\subsection{Other techniques to reduce memory consumption}

A more complex simplification comes from noticing that a signature of the form
$a[b]c$ will not touch any of the $b$ values until all of $c$ is dealt with. This
means that there will be a set of signature prefixes $p_i$ with multiplicities $m_i,$ independent of $b,$ such that
$f(a[b]c)$ = $\sum\limits_{i} m_i f(p_i b)$. These $p_i$ and $m_i$ can be computed when needed and
cached. If you did this for all $a[b]c$ you would spend more time and memory
on this optimization than the original problem involved. However, if you just
do it for sufficiently short $a$ (in practice we used length of $a < 8$) then it can
reduce memory consumption by about 30 percent without significant effect on speed.
Note that the length of $a$ is more important
than the length of $c$, as all prefixes $p_i$ must be no longer than $a$, so a small length
of $a$ ensures that the number of terms here does not get too large. Also note that
$[b]$ can actually be a series of brackets, e.g. $[b_1][b_2][b_3][b_4]$.

Reducing the number of signatures reduces both the memory and the execution time, but
memory consumption tends to be the bottleneck. One simple trick for problems of this
class is to only save the result some fraction $p$ of the time. Signatures that are
only used once may not take up memory, and frequently used signatures will get
stored eventually. The smaller $p$, the less memory used, but the more time used. We
found that $p=0.3$ reduced memory consumption by about 30 percent with about a 30 percent
increase in execution time. This indicates that a significant fraction of the values $f(S)$
computed are only used once.

\subsection{Implementation}

The algorithm was implemented in Scala, which compiles to Java virtual machine bytecode,
and run on a computer with 1TB RAM.

The signature described here is less straightforward to encode on a computer than
the array of 32 bytes in \cite{JN13}. However, in practice it can be easily done in 128 bits
for all the signatures needed for $n$ up to the 50s.
This makes the keys 16 bytes, reducing memory consumption.

The signatures are encoded as a bitstring. The bitstring starts with a 6 bit number
equal to the sum of the integers in the signature $(n).$ This is used for determining
when the bitstring stops. There is then one bit to indicate whether the signature starts
with a comma. Then there is a repeating series of 1 bit for whether there is an
open bracket, 2 to 8 bits for encoding an integer, and one bit
for whether there is a closing bracket after that integer. This repeats until everything
is encoded.

Integers were encoded as follows:

\begin{itemize}

  \item 00 means the number 1,

  \item 01 means the number 2,

  \item 10bbb means the number bbb+3,

  \item 11bbbbbb means the number bbbbbb+11

\end{itemize}

Some signatures have a large number of small numbers in them; this will encode them in
a small number of bits. Others have a small number of large numbers; this will also encode
them in a small number of bits. Indeed, up until $n=29$, only 64 bits are needed for
the keys.

The values are stored as 64 bit values {\em modulo} some large prime. The computation is redone
{\em modulo} a different prime and the values reassembled using the Chinese remainder theorem.
This somewhat reduces memory use relative to 128 bit values, but also simplifies implementation
as the java virtual machine does not handle 128 bit integers easily.

The main memoization store was therefore effectively a large 128 bit key to 64 bit value
hash table. As the upper 64 bits of the key were sparsely used (indeed only 131 different
values for 35 terms), a significant memory saving was generated by having a master hash map
from the upper 64 bits to a slave hash map. The slave hash map then only needed to use the
lower 64 bits as keys. In practice the slave hash maps were still further subdivided in
one more layer; this was due to java virtual machine limitations on array size restricting
the length of a hash map. But the end result is the main memory store used 64 bit keys
and values in a gnu trove hash map.

Various of the tricks used increased computation time to a matter of days; we made a
parallelized version which worked about twice as fast as the single threaded version
on a 4 core desktop pc for testing, but turned out slower on the 32 core production
machine than the single threaded version, possibly because of memory coherency overheads.
So it was not used.
The parallelization comes from having different threads compute different signatures.
If, in computing $f(S)$, the value $f(s)$ for a subsignature $s$ is not available, then
a placeholder is inserted in $f(S)$, and $s$ is added to a priority queue of signatures
to evaluate. When a processor has nothing else to do, it takes the signature with the
lowest $n$ from the priority queue and evaluates it. When a signature is finished, it
fills in all the placeholders for that signature. Taking the signature with the
lowest $n$ prevents the queue from growing exponentially with $n$.

 With all these improvements implemented, the algorithm produced 5 further terms on a comparable computer. The calculation had to be run twice, {\em modulo} two different primes, and the result reconstructed by the Chinese Remainder Theorem. As a check, we also ran it a third time, {\em modulo} another prime.
For each of the two moduli used, the program ran for 5 days and used somewhat over 540 GB of memory. This is the amount of memory used at the termination of the algorithm as reported by the JVM; some extra is needed for various tasks, primarily temporary objects such as resizing of hash tables.

Source code is available on https://github.com/AndrewConway/enumeration/ in the folder
avoid1324.

The coefficients are given in Table \ref{tab:coeffs}.

 \begin{table}[htbp] 
   \centering
   \caption{Coefficients of $1324$ pattern-avoiding permutations.}
\label{tab:coeffs}
\begin{center}
\begin{tabular}{ll} \hline \hline
   1& 3421888118907\\
1&25887131596018\\
2&198244731603623\\
6&1535346218316422\\
23&12015325816028313\\
103&94944352095728825\\
513&757046484552152932\\
2762&6087537591051072864\\
15793&49339914891701589053\\
94776&402890652358573525928\\
591950&3313004165660965754922\\
3824112&27424185239545986820514\\
25431452&228437994561962363104048\\
173453058&1914189093351633702834757\\
1209639642&16130725510342551986540152\\
8604450011&136664757387536091240503406\\
62300851632&1163812341034817216384582333\\
458374397312&9959364766841851088593974979\\
&85626551244475524038311935717\\
  \hline
   \end{tabular}
   \end{center}
 \end{table}

\section{Analysis}
In the case of a simple algebraic singularity with asymptotic form (\ref{eq:alg}), the ratio of the coefficients is
\BE \label{eq:rat1}
r_n = \frac{a_n}{a_{n-1}} \sim \mu\left (1 + \frac{g}{n} +{\rm O}(\frac{1}{n^2}) \right ).
\EE
If on the other hand the coefficients of some generating function are as in eqn. (\ref{eq:can1}), then the ratio of successive coefficients $r_n = b_n/b_{n-1},$ is
\begin{multline} \label{eq:rn}
r_n = \mu \left (1 + \frac{\sigma \log \mu_1}{n^{1-\sigma}} + \frac{g}{n} + \frac{\sigma^2 \log^2 \mu_1}{2n^{2-2\sigma}} + \frac {(\sigma-\sigma^2)\log \mu_1+2g\sigma \log \mu_1}{2n^{2-\sigma}} \right . \\
 \left . {}+ \frac{\sigma^3 \log^3 \mu_1}{6n^{3-3\sigma}} +{\rm O}(n^{2\sigma-3}) + {\rm O}(n^{-2}) \right ).
\end{multline}
In particular,
when $\sigma = \frac{1}{2},$ this specialises to
\BE \label{eq:half}
r_n = \mu \left (1 + \frac{ \log \mu_1}{2\sqrt{n}} + \frac{g+\frac{1}{8}\log^2 \mu_1}{n} + \frac{\log^3\mu_1+(6+24g)\log \mu_1 }{48n^{3/2} } + {\rm O}(n^{-2}) \right ).
\EE

 \subsection{Differential approximant analysis}
The most successful numerical method for extracting the asymptotics from the first few terms of the OGF of a function with an algebraic singularity\footnote{With a slight abuse of notation, we refer to a singularity of the form $(1-x/x_c)^{\alpha}$ as an algebraic singularity even in those cases where $\alpha$ is not rational.} is the method of differential approximants, (called the DA method) due to Guttmann and Joyce \cite{GJ72}, with subsequent refinements due to Baker and Hunter \cite{HB79} and Fisher and Au-Yang \cite{FA79}. Details are given in \cite{G14, GJ09, G89}. In brief, the method fits available coefficients to a judiciously chosen family of D-finite ordinary differential equations (ODEs), and the singularity structure of the ODEs is extracted by standard methods \cite{Ince27, Forsyth02}. 

For models with an isolated algebraic singularity, the method is very successful, with the radius of convergence and critical exponents frequently estimated to 10 significant digit accuracy or better, from a series of length 30-80 terms. However, when the method is used to analyse models with singularities that are not algebraic, such as those whose coefficients have the asymptotic form (\ref{eq:can1}), the method fails, though in a predictable manner. That is to say, one finds that the radius of convergence estimates are typically only found to two or three significant digits, and the critical exponent estimates are numerically large, typically around 10 or -10. 

In this way, the method is useful -- as is a canary in a coal mine. If one analyses the known terms of the series with the method of differential approximants and finds estimates of the radius of convergence to be poorly converged, with numerically large exponent values, one can be confident that the underlying OGF does not have an algebraic singularity. Applying the method to the first 30 terms of the 1342 and 1234 PAPs, the known solutions are found. Applying the method to the 36 coefficients we have for 1324 PAPs, the method suggests that the radius of convergence is around 0.09, with an exponent variously estimated to be -20 or +15 or anything in between! This is the hallmark of a non-algebraic singularity. 

Further details of the DA method, its successes and limitations are discussed in \cite{G14}. For the moment, we simply conclude that the OGF of 1324 PAPs almost certainly does not have an algebraic singularity. In the next section we explore the nature of the singularity by looking at the ratio of successive coefficients. 

\subsection{Ratio analysis}
\begin{figure}[t!]
\centering
\subfigure[Plot of ratios of coefficients against $\frac{1}{n}$.]{\includegraphics[width=7.3cm]{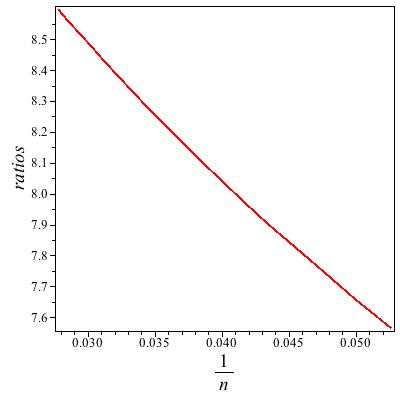}} 
\subfigure[Plot of ratios of coefficients  against $\frac{1}{\sqrt{n}}$.]{\includegraphics[width=7.3cm]{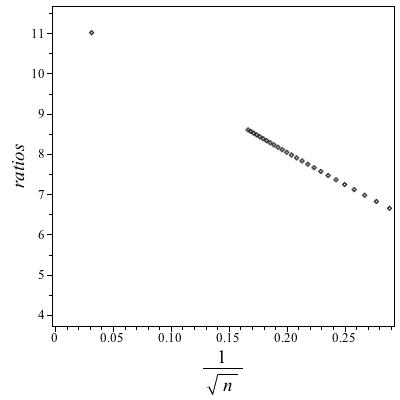}} 
\caption{}
\label{fig:ratios}
\end{figure}

In order to determine the nature of the asymptotic form of the coefficients of the 1324-PAP OGF, 
we first plot the ratios of successive coefficients $r_n = p_n/p_{n-1}$ against $1/n,$ as shown in figure \ref{fig:ratios}(a).  The locus is clearly concave. This is inconsistent with an algebraic singularity, as can be seen from eqn. (\ref{eq:rat1}). We next plot the same ratios against $1/\sqrt{n}$ in figure \ref{fig:ratios}(b), this time as a point plot, and the plot is seen to be visually linear, implying, from eqn (\ref{eq:rn}) that $\sigma \approx 1/2.$ The outlying point near to the vertical axis is a Monte Carlo result obtained by Steingr\'imsson \cite{ES12} for PAPs of length 1001.

\begin{figure}[t!]
\centering
\subfigure[Plot  of extrapolated ratios against $\frac{1}{n}$.]{\includegraphics[width=7.3cm]{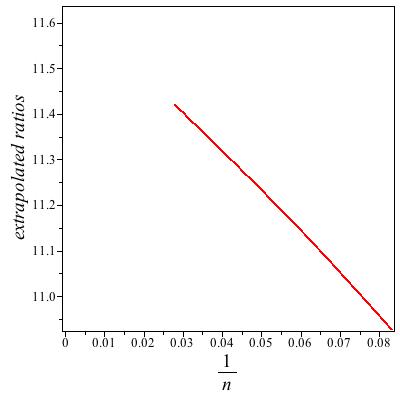}} 
\subfigure[Plot of estimators of $\mu$ by square root intercepts.]{\includegraphics[width=7.3cm]{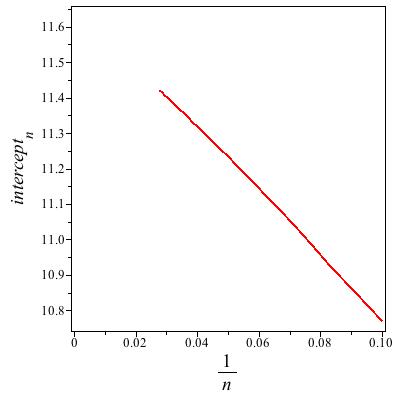}} 
\caption{}
\label{fig:sig1}
\end{figure}

 Linear extrapolation (not including the isolated point) implies a limiting value as $n \to \infty$ around 11.5. We can significantly improve on this estimate by considering the sequence of extrapolants defined by successive pairs of points. That is to say, one can simply linearly extrapolate successive pairs of ratios $(r_k,r_{k+1})$ with $k$ increasing up to the maximum value achievable with our data, which is 35. A plot of successive extrapolants against $1/n$ is shown in figure \ref{fig:sig1}(a), which appears to be linear. A crude extrapolation with a ruler suggests a limit of around 11.60. 

Assuming (tentatively) that $\sigma = 1/2,$  another way to use the ratios to get a better estimate of $\mu$ is to eliminate the assumed term O$(1/\sqrt{n})$ in the ratios by forming the modified ratios
\BE \label{eq:modrat}
intercept_n = \frac{\sqrt{n} \cdot r_n - \sqrt{n-1} \cdot r_{n-1}}{\sqrt{n} -\sqrt{n-1}} \sim \mu \cdot \left ( 1 + {\rm O} (1/n) \right).
\EE
We show in Figure \ref{fig:sig1}(b) a plot of $intercept_n$ against $1/n,$ and this too appears to be going to a value very close to 11.60.
 We will take that as our initial estimate, which we subsequently refine. We also take $\sigma =1/2$ as our (initial) conjectured value. In doing so we are, in part, relying on the observation that in all known cases \cite{G14} where this asymptotic behaviour is observed, $\sigma$ is a simple rational, usually $1/2$ or $1/3$.

We also plotted (but don't display) the ratios against $1/n^{2/3},$ which is appropriate if $\sigma = 1/3.$ In that case the locus was convex, rather than concave. So on the basis of ratio plots alone, $\sigma \approx 1/2$ is our tentative estimate. 

\begin{figure}[t!]
\centering
\subfigure[Log-log plot  of $(1-r_n/\mu)$ against $n$..]{\includegraphics[width=7.3cm]{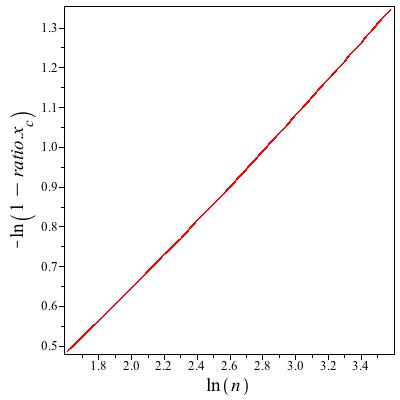}} 
\subfigure[Plot of local extrapolants of figure at left, estimating $1-\sigma.$]{\includegraphics[width=7.3cm]{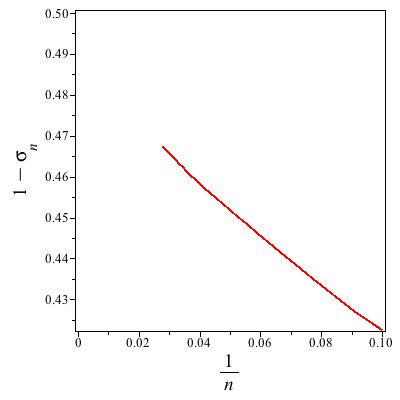}} 
\caption{}
\label{fig:sig2}
\end{figure}

In order to more accurately determine the value of the exponent $\sigma$, we note from (\ref{eq:rn}) that $$(r_n/\mu-1) \sim const. n^{\sigma-1}.$$ We show in figure \ref{fig:sig2}(a) a log-log plot of $(1-r_n/\mu)$ against ${n},$ where we have taken 11.60 as the (tentative) value of $\mu.$ This should be linear, with gradient $\sigma-1.$ A small degree of curvature is evident. Accordingly, we extrapolate the local ratios, defined as
\BE \label{eq:locrat}
1-\sigma_n = \frac{\log \left ( 1 -\frac{r_{n-1}}{\mu} \right )-\log \left ( 1 -\frac{r_{n}}{\mu} \right )}{\log {n} -\log(n-1)},
\EE
against $1/n.$ The results are shown in Figure \ref{fig:sig2}(b)
The ordinates are estimators of $1-\sigma.$ If one accepts that $\sigma$ is  a simple rational number, the value $1/2$ is inescapable. 

\begin{figure}[t!]
\centering
\subfigure[Log-log plot  of $r_{\sigma_n}$ against $n$.]{\includegraphics[width=7.3cm]{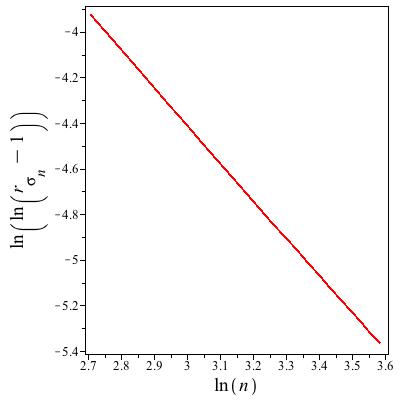}} 
\subfigure[Local gradient of log-log plots  of  $r_{\sigma_n}$  against $\frac{1}{n}$.]{\includegraphics[width=7.3cm]{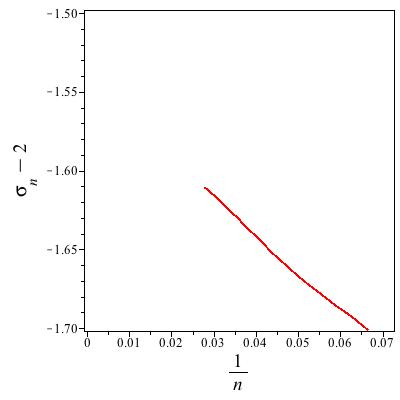}} 
\caption{}
\label{fig:sig3}
\end{figure}

We can also estimate $\sigma$ without assuming the value of $\mu$ as follows. From eqn. (\ref{eq:can1}), one sees that
\BE \label{eq:sig1}
r_{\sigma_n} = \frac{b_n \cdot b_{n-2}}{b_{n-1}^2} \sim 1 + \frac{\sigma \cdot (\sigma-1)\log{\mu_1}}{n^{2-\sigma}} + {\rm O}(1/n^2),
\EE
so $\sigma$ can be estimated from a log-log plot of $r_{\sigma_n}-1$ against $n,$ independent of the value of $\mu.$ From this plot, shown in Figure \ref{fig:sig3}(a), which is linear as expected, we calculate the local gradient at each pair of successive points, as described above, and plot these against $1/n.$
This  is shown in figure \ref{fig:sig3}(b), and it can be seen that the expected limit, as $n \to 0$ is, plausibly, $-1.5.$ From eqn. (\ref{eq:sig1}), this limit should be $2 - \sigma,$  which is consistent with our assertion that $\sigma = 1/2.$
In our subsequent analysis, we will assume this value.

With $\sigma$ taken to be $1/2,$ we attempt to refine the estimate  of $\mu$ by extrapolating the ratios of the coefficients $r_n$ using the Bulirsch-Stoer \cite{BS64} algorithm, with parameter $w=1/2.$ This algorithm extrapolates a sequence $\{s_n\}$ assuming $s_n \sim s_\infty + c/n^w + o(n^{-w}),$ where $w$ is provided by the user. The results are given in Table \ref{tab:BS3}. Each successive row represents a higher order of extrapolation. We only show the last 7 entries of each order of extrapolation. We continue generating rows until the row entries lose monotonicity. From the first row we conclude $\mu < 12.832.$ From the second we conclude $\mu < 11.6663,$ and from the third we conclude $\mu < 11.6112.$  The rate of decrease in the entries in the third row is consistent with our initial estimate of 11.60, so we will retain this estimate for the time being, as the extrapolation has not given us a more precise value. It does however add support to that choice.

\begin{table}
\caption{\label{tab:BS3}
Last seven entries in each row of the table of Bulirsch-Stoer extrapolants. Each successive row is the result of a successively higher degree of extrapolation. The ratios of successive coefficients are extrapolated here, with parameter $w=1/2.$}
\begin{center}
\begin{tabular}{llllllll} \hline \hline
 L   &  T(L,N-L-6) & T(L,N-L-5) & T(L,N-L-4)& T(L,N-L-3) & T(L,N-L-2) & T(L,N-L-1) & T(L,N-L) \\
 
\hline

\hline
1&13.17792074& 13.10768870 &13.04333009 &12.98413360 &12.92949766 &12.87891071& 12.83193533 \\
2&11.67944917 &11.67746701 &11.67538180 &11.67320098 &11.67094284 &11.66863183 &11.66629498 \\
3&11.64370166 &11.63846233 &11.63275886 &11.62693780 &11.62127955 &11.61599904 &11.61124940 \\
4&11.73026868 &11.71850728 &11.71746270 &11.72167941 &11.73128989 &11.74947827 &11.78630320\\

\hline \hline
\end{tabular}
\end{center}
\end{table}

Assuming then that $\sigma=1/2,$ from (\ref{eq:half}), it follows that $$r_n/\mu = 1 + \frac{ \log \mu_1}{2\sqrt{n}} + \frac{g+\frac{1}{8}\log^2 \mu_1}{n} + {\rm O}(n^{-3/2}). $$

In order to estimate $\mu_1$ and $g,$ we solve, sequentially, the trio of equations 
\BE \label{eq:fit2a}
r_j/\mu = 1 + \frac{ c_1}{\sqrt{j}} + \frac{c_2}{j} +  \frac{c_3}{j^{3/2}}, 
\EE
for $ j=k-1,$ $j=k$ and $j=k+1,$ with $k$ ranging from 2 up to 35, and $\mu$ set at 11.60.

\begin{figure}[t!]
\centering
\subfigure[Plot of  parameter $c_1$ of (\ref{eq:fit2a}) against $\frac{1}{n^{3/2}}$.]{\includegraphics[width=7.3cm]{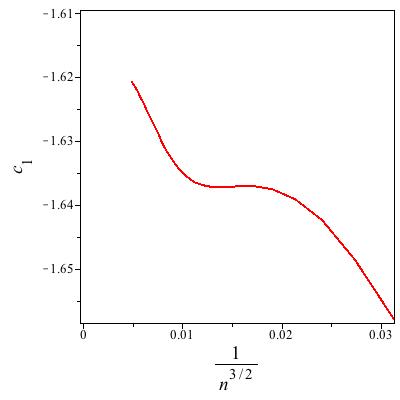}} 
\subfigure[Plot of  parameter $c_2$ of (\ref{eq:fit2a}) against $\frac{1}{n}$.]{\includegraphics[width=7.3cm]{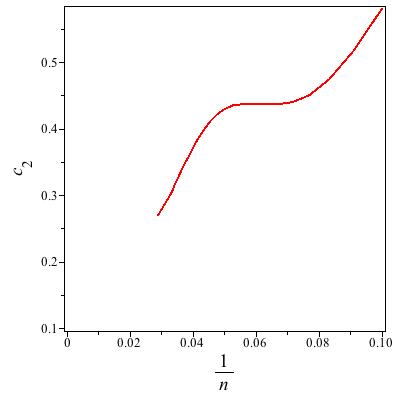}} 
\caption{}
\label{fig:papc}
\end{figure}

The results are shown in figures \ref{fig:papc}(a) and \ref{fig:papc}(b), plotting the parameters $c_1$ and $c_2$ respectively. The first neglected term in the asymptotics is O$(1/n^2)$ which is O$(1/n^{3/2})$ smaller than the term with coefficient $c_1,$ so $c_1$ is plotted against $1/n^{3/2}.$ By a similar argument, $c_2$ is plotted against $1/{n}.$ A simple visual extrapolation gives the estimate $c_1 = -1.615 \pm 0.005.$ The plot for $c_2$ is difficult to extrapolate. It appears to be turning near its end point, and we very tentatively estimate  $c_2 \approx 0.15.$ Unless the gradient changes sign, we can only say that $c_2 < 0.2,$ and it seems to be going to a positive value. If the gradient changes sign, we can't even say that. We don't show the plot of $c_3$ as we cannot extrapolate it. From (\ref{eq:half}), $c_1=\log{\mu_1}/2$ and $c_2=g+\frac{1}{8}\log^2 \mu_1.$ Hence we estimate  $\log{\mu_1} \approx -3.23,$ and assuming $c_2$ is in the range $[0,0.2]$ this gives $g = -1.2 \pm 0.15.$ We repeated this analysis varying the estimate of $\mu$ in the range [10.58,11.62]. With $\mu$ in this range, we estimate $\log{\mu_1} = 3.23 \pm 0.07,$ and $g = -1.2 \pm 0.4.$

An alternative form of analysis involves direct fitting to the parameters in the assumed asymptotic form. That is to say, the assumed asymptotic form is $$b_n \sim B\cdot \mu^n \cdot \mu_1^{n^\sigma} \cdot n^g.$$ Therefore
\BE \label{logcan1}
\log {b_n} \sim \log{B} + n \log{\mu} + n^\sigma \log{\mu_1} + g \log {n}.
\EE
So if $\sigma$ is known, or assumed, we have four unknowns in this linear equation. It is then straightforward to solve the linear system
$$\log {b_k} =  c_1k  + c_2 k^\sigma  +c_3 \log {k}+c_4$$ for $k=n-2,\, n-1, \, n, \, n+1$ with $n$ ranging from $3$ to $35.$ Then $c_1$ estimates $\log {\mu},$ $c_2$ estimates $\log{\mu_1}$, $c_3$ estimates $g$ and $c_4$ gives estimators of $\log{B}.$ 

An obvious useful variation is in those cases where, say, $\mu$ is known, or accurately estimated. Then one can solve  
\BE \label{eq:three}
\log b_n - \mu \log{n} =  c_1 n^\sigma  +c_2 \log {n} +c_3
\EE
 from three successive coefficients $b_{n-1}, \,\, b_n, \,\, b_{n+1}$, as before increasing the order of the lowest used coefficient by one until one runs out of coefficients. We do this below with $\mu$ varying within its estimated error range.

Fitting the available coefficients to the four unknowns, we estimate $c_1\approx 2.450 \pm 0.002,$ implying $\mu= 11.59 \pm 0.02,$ (in good agreement with our earlier estimate of 11.60), $c_2 =-3.23 \pm 0.03,$ implying $\mu_1 = 0.0396 \pm 0.0012,$   $c_3\approx -1,$ while $c_4$ is difficult to estimate beyond saying it is in the range $[1.3,3]$ implying $B \in [4,20].$ 

We repeated this analysis with a 3 parameter fit, varying $\mu$ in the range $[11.58,11.62].$ This gave $c_2 =-3.22 \pm 0.08,$ implying $\mu_1 = 0.040 \pm 0.003,$   $c_3=g= -1.15 \pm 0.2,$ and $c_4=1.8 \pm 0.5$ implying $B = 7 \pm 3.$ 

As noted above, differential approximants are useful insofar as they indicate that the singularity is not algebraic. They provide a signal, but are then of no further use in their current form. The presence of the $\mu_1^{\sqrt{n}}$ term is responsible for the lack of applicability of the method. However we can manipulate the series to remove the offending term, and then use this powerful method.  From eqn. (\ref{logcan1}), defining $\tilde{b}_n = b_n/\sqrt{n},$ one has

\BE  \label{eq:renorm}
c_n=2n^{3/2}(\tilde{b}_n - \tilde{b}_{n-1}) \sim (2g-\log{B}) + n\log(\mu)-g\log(n)
\EE
So 
$$d_n = \exp(c_n) \sim D\cdot \mu^n \cdot n^{-g} \cdot (1+ {\rm O}(n^{1-\sigma})),$$ with $D=\exp(2g)/B.$

The dominant asymptotics of the coefficients $d_n$ are now those of an algebraic singularity -- though the nature of the correction term, O$(n^{1-\sigma}),$ means that there is a confluent singularity with exponent less than 1. This means that the OGF $\sum d_n\cdot x^n$ can be analysed by standard methods used to analyse series with algebraic singularities. This includes the method of differential approximants and Bulirsch-Stoer extrapolation of ratios. For the former method, one should use 3rd order DAs, as the presence of a confluent singularity means that we require an ODE with two independent solutions, and to allow for a non-singular background term requires a third independent solution. For Bulirsch-Stoer extrapolation of ratios, the parameter $w=1$ should be used, as the ratios $d_n/d_{n-1} \sim \mu \cdot (1 -g/n + {\rm o}(1/n)).$ 

 For the DA analysis we have used 3rd and 4th order ODEs. As an aside, we also tried 2nd order ODEs and these were unsatisfactory, as we expected, as they produced two singularities close together in an unsuccessful attempt at representing the confluent singularity. We summarise the results in table \ref{tab:ana}. The column labelled $L$ gives the degree of the inhomogeneous polynomial of the approximating ODEs. The entries give estimates, averaged over many approximants, of the position and exponent of the singularity of the  ODEs. Full details of the method are given in \cite{G89, GJ09}. It is seen that the 3rd order DAs give estimates of the radius of convergence that are centred around $0.086140 = 1/11.609,$ with exponent estimates around $g \approx -0.93.$  The 4th order approximants give slightly higher estimates of both the critical point and the absolute value of the exponent. We estimate $
1/\mu \approx 0.08619,$ or $\mu \approx 11.602.$ This is remarkably close to the initial estimate above, $\mu =11.60.$ For the exponent $g$ we estimate $g = -1.0 \pm 0.1.$

\begin{table}
\caption{\label{tab:ana}
Critical point and exponent estimates for renormalised 1324 PAPs}
\begin{center}
\begin{tabular}{lllll} \hline \hline
 $L$   &  \multicolumn{2}{c}{Second order DA} &
       \multicolumn{2}{c}{Third order DA} \\
\hline
    &  \multicolumn{1}{c}{$1/\mu$} & \multicolumn{1}{c}{$g-1$} &
      \multicolumn{1}{c}{$1/\mu$} & \multicolumn{1}{c}{$g-1$} \\
\hline
 0 & 0.086237& -2.160 &  0.086134& -1.925 \\
 1 & 0.086116& -1.905 &  0.086156& -1.958 \\
2 & 0.086142& -1.942 & 0.086149& -1.944 \\
3 & 0.086159& -1.964 & 0.086167& -1.982 \\
4 & 0.086137& -1.926 & 0.086156& -1.960 \\
 5 & 0.086111& -1.914 &  0.086162& -1.966 \\
 6& 0.086110& -1.910 &  0.086158& -1.959 \\
7 & 0.086140& -1.926 & 0.086170& -1.980 \\
8 & 0.086143& -1.936& 0.086178& -1.997 \\
9 & 0.086145& -1.934 & 0.086186& -2.012 \\
10 & 0.086142& -1.931 & 0.086188& -2.011 \\
\hline \hline
\end{tabular}
\end{center}
\end{table}

We can also apply other standard techniques to the transformed series. The ratios of successive terms ($d_n$) of the transformed series when plotted against $1/n$ are now visually linear. Accordingly, we extrapolate the ratios of the coefficients of the transformed series using the Bulirsch-Stoer algorithm, with parameter $w=1.$ The results are shown in table \ref{tab:BS4}. The first two rows are behaving monotonically. The last entry in the second row suggests that $\mu < 11.622,$ and assessing the way entries in that row are decreasing, we judge the limit to be around $11.60 \pm 0.01.$ Together with the result of the DA analysis given above, we combine these two results and give as our final estimate $\mu=11.60 \pm 0.01.$ The quoted error is to be interpreted as a confidence interval, not a rigorous error bound. And, it should be stressed, our analysis is predicated on our assumption that $\sigma =1/2.$

\begin{table}
\caption{\label{tab:BS4}
Last seven entries in each row of the table of Bulirsch-Stoer extrapolants. Each successive row is the result of a successively higher degree of extrapolation. The ratios of successive coefficients of the transformed series are extrapolated here.}
\begin{center}
\begin{tabular}{llllllll} \hline \hline
 L   &  T(L,N-L-6) & T(L,N-L-5) & T(L,N-L-4)& T(L,N-L-3) & T(L,N-L-2) & T(L,N-L-1) & T(L,N-L) \\
 
\hline

\hline
1&11.63526435& 11.63607983 &11.63648194& 11.63651932 &11.63624657& 11.63571995 &11.63499412\\ 
2&11.65200739 &11.64744495 &11.64238782 &11.63709746 &11.63181012 &11.62672996 &11.62202312 \\
3&11.64102224& 11.63870198& 11.63726970 &11.63654375& 11.63647001 &11.63709518& 11.63859791\\

\hline \hline
\end{tabular}
\end{center}
\end{table}

Finally, we estimate the amplitude $B$ by extrapolating the sequence $b_n/(\mu^n \cdot \mu_1^{\sqrt{n}} \cdot n^g),$ supplying the estimates of the critical parameters already found, against $1/n.$ In this way we estimate $B = 9.5 \pm 0.5.$

\section{Conclusion}
We have given a refined version of the JN algorithm that allows five further coefficients of the 1324 PAP generating function to be obtained, with comparable computer resources. 

Analysing the coefficients of the generating function, we provide compelling evidence that they have a singularity structure of the form
$$B\cdot \mu^n \cdot \mu_1^{\sqrt{n}} \cdot n^g.$$ We give as our final estimates of the critical parameters $\mu=11.60 \pm 0.01,$ $\sigma=1/2,$ $\mu_1 = 0.0398 \pm 0.0010,$ $g = -1.1 \pm 0.2$ and $B =9.5 \pm 1.0.$ If, as is seen in other problems which have coefficients with a similar asymptotic structure, that $g$ is a simple rational fraction, the most likely is $-7/6$ or $-6/5,$ though we could not rule out $-1$ or even $-5/4.$

 Zeilberger has said ``Not even God knows the number of 1324-avoiders of length 1000". While making no Messianic claims, our asymptotics permit the approximate answer $4.6 \times 10^{1017}.$

\section*{Acknowledgements}
AJG wishes to acknowledge helpful conversations with Einar Steingr\'imsson, who brought this problem to our attention, and the hospitality of Mathematisches Forschungsinstitut Oberwolfach  and the Enumerative Combinatorics Workshop held there on March 2--8, 2014, where this work was initiated. We are grateful to Alan Sokal, who gave us access to his Dell computers with 1TB memory, which were needed for these calculations, through his NSF grant PHY--0424082 NYU.    AJG wishes to thank the Australian Research Council for supporting this work through grant DP120100931.

\end{document}